\documentclass[11pt]{article}
\usepackage{graphicx,epsfig}
\usepackage{amsmath,amssymb,euscript}
\usepackage{latexsym}
\usepackage{color}
\newtheorem{theorem}{\bf Theorem}[section]

\newtheorem{lemma}[theorem]{\bf Lemma}

\newcommand{\Proofend}{\hfill$\Box$\vspace{.4cm}}
\newcommand{\interior}[1]{{\kern0pt#1}^{\mathrm{o}}}

\def\NN{{\mathbb N}}
\def\RR{{\mathbb R}}

\def\PP{{\mathbb P}}
\def\SS{{\mathbb S}}
\def\C{{\mathcal{C}}}

\def\V{{\mathcal{V}}}

\def\S{{\mathcal{S}}}

\def\V{{\mathcal{V}}}
\def\x{{\bf x}}
\def\y{{\bf y}}

\DeclareMathOperator{\sign}{sign}

\newcommand{\cmarc}[1]{{\color{black} #1}}

\begin{document}
\title{Uniform approximation on the sphere by least squares polynomials}
\author{Woula Themistoclakis\thanks{C.N.R. National
        Research Council of Italy,
        Istituto per le Applicazioni del Calcolo ``Mauro Picone'',  via P. Castellino, 111, 80131 Napoli, Italy.
        {\tt woula.themistoclakis@cnr.it}.}
        \and Marc Van Barel\thanks{KU Leuven, Department of Computer Science, KU Leuven,
Celestijnenlaan 200A,
B-3001 Leuven (Heverlee), Belgium. {\tt marc.vanbarel@cs.kuleuven.be}.
}}
%\date{}
\maketitle
\begin{abstract}
The paper concerns the uniform polynomial approximation of a function $f$, continuous on the unit Euclidean sphere of $\RR^3$ and known only at a finite number of points that are somehow uniformly distributed on the sphere.
First we focus on least squares polynomial approximation and prove that the related Lebesgue constants w.r.t.\ the uniform norm grow at the optimal rate. Then, we consider delayed arithmetic means of least squares polynomials whose degrees vary from $n-m$ up to $n+m$, being $m=\lfloor \theta n\rfloor$ for any fixed parameter $0<\theta<1$. As $n$ tends to infinity, we prove that these polynomials uniformly converge to $f$ at the near-best polynomial approximation rate. Moreover, for fixed $n$, by using the same data points we can further improve the approximation by suitably modulating the action ray $m$ determined by the parameter $\theta$. Some numerical experiments are given to illustrate the theoretical results.
\\[0.2cm]
{\bf keywords:} polynomial approximation on the sphere, least squares approximation, uniform approximation, Lebesgue constant,
de la Vall\'ee Poussin type mean.
\\[0.2cm]
\end{abstract}
\section{Introduction}
Many applications in various fields of physics, biology and engineering, require a fast polynomial reconstruction of a real-valued\footnote{The generalization to complex-valued functions is straightforward.} function $f$ defined on the sphere $\SS^2=\{\x:=(x,y,z)\in\RR^3: x^2+y^2+z^2=1\}$, by using sampled values of $f$ at a discrete point set
\[
X_N:=\{\xi_1,\ldots,\xi_N\}\subset\SS^2.
\]
%In this paper we focus on polynomial approximation which turns to be useful in several practical applications, being polynomials the eigenfunctions of many pseudodifferential operators, \cmarc{further than infinitely smooth.}

We assume that the point set $X_N$ is distributed on the sphere in such a way to support a positive weighted quadrature rule of a suitable degree of exactness $\mu\in\NN$, i.e. there exist positive real numbers $\lambda_1,\ldots,\lambda_N$ such that
\begin{equation}\label{quad}
\int_{\SS^2} f(\x)d\sigma(\x)=\sum_{i=1}^N\lambda_i f(\xi_i),\qquad \forall f\in\PP_{\mu},
\end{equation}
where $d\sigma$ denotes the usual surface measure on $\SS^2$ and $\PP_\mu$ is the space of all spherical polynomials (i.e., polynomials in three variables restricted to $\SS^2$) of degree at most $\mu$.

Moreover, we suppose that the following Marcinkiewicz type inequality holds
\[
\frac 1{\mu^2}\sum_{i=1}^N|f(\xi_i)|\le \C \int_{\SS^2}|f(\x)|d\sigma(\x),\qquad \forall f\in\PP_\mu,\qquad \C\ne\C(f,\mu,N),
\]
where throughout the paper, $\C$ denotes a positive constant, which can take on different values in the different positions where it appears, and we write $\C\ne\C(f,\mu,...)$  to mean that $\C$ is independent of $f,\mu,...$.

Examples as well as sufficient conditions on $X_N$ for the existence of positive quadrature rules have been derived by many authors (see, e.g., \cite{r943,b509,r938,r939,W2,r947}).

In literature \cite{W4}, the positive weighted quadrature rule (\ref{quad}) with $\mu=2n$ has been applied to the Fourier orthogonal projection $\S_n:f\rightarrow \S_nf\in\PP_n$ w.r.t. the scalar product
\begin{equation}\label{pr-cont}
<f,g> :=\int_{\SS^2}f(\x)g(\x)d\sigma(\x),
\end{equation}
in order to get the so--called {\it hyperinterpolation polynomial}, which we denote by $L_nf$.

Also, positive weighted quadrature rules have been used in order to discretize some generalized de la Vall\'ee Poussin means \cite{r938, r947} (i.e., delayed weighted means of Fourier sums), obtaining the so--called {\it filtered hyperinterpolation polynomials} \cite{r951}.

Concerning the  approximation properties,  both hyperinterpolation and filtered hyperinterpolation polynomials are comparable with their respective continuous versions.

More precisely, it is known \cite{W-Da} that Lebesgue constants of Fourier partial sums $S_nf$ grow with the degree $n$ at the minimal projection rate, that is $\sqrt{n}$, and similarly, for the hyperinterpolation polynomial $L_nf$ we have \cite{W-Re, W4}
\begin{equation}\label{LC-Fou}
 \|L_n\|\sim\|\S_n\|\sim \sqrt{n},
\end{equation}
where we set $\|T\|:=\sup_{\|f\|_\infty\le 1}\|Tf\|_\infty$ with $\|f\|_\infty:=\sup_{\x\in\SS^2}|f(\x)|$, and by $a_n\sim b_n$ we mean that $c_1 a_n\le b_n\le c_2 a_n$ with $c_1,c_2>0$ independent of $n$.

Moreover, it is known that we can get uniformly bounded Lebesgue constants by filtered hyperinterpolation quasi--projections, under suitable assumptions on the filter coefficients which define the generalized de la Vall\'ee Poussin mean \cite{r938, r951, W1}.\newline
In particular, for any $0<\theta <1$, if we set $m=\lfloor \theta n\rfloor$ ($\lfloor \cdot\rfloor$ being the floor function) and apply (\ref{quad}) with $\mu=4n$ to the following arithmetic mean of Fourier sums
\begin{equation}\label{VP-cont}
\V_n^mf(\x)=\frac 1{2m+1}\sum_{r=n-m}^{n+m}\S_rf(\x),\qquad m=\lfloor \theta n\rfloor, \qquad \x\in\SS^2,
\end{equation}
then the resulting filtered hyperinterpolation polynomial, which we denote by $V_n^mf$, satisfies \cite{W1}
\begin{equation}\label{LC-VP}
\|V_n^m\|\le\C \|\V_n^m\|\le \C,\qquad n\in\NN,\qquad m=\lfloor \theta n\rfloor,\qquad \C\ne\C(n).
\end{equation}
The estimates (\ref{LC-Fou}) and (\ref{LC-VP}) imply that, as $n\rightarrow\infty$, filtered hyperinterpolation polynomials $V_n^mf$ uniformly converge to $f$ for any continuous function $f$ on the sphere, while hyperinterpolation polynomials $L_nf$ do it if the function $f$ is smooth enough to satisfy $\lim_{n\rightarrow\infty}\sqrt{n}E_n(f)=0$, with
\begin{equation}\label{En}
E_n(f):=\inf_{P\in\PP_n}\|f-P\|_\infty,
\end{equation}
the error of best uniform approximation of $f$ in $\PP_n$.

On the other hand, we observe that both the discrete approximation polynomials  $L_nf$ and $V_n^mf$ have been obtained by applying a suitable positive quadrature rule and, for this reason, they are based not only on the function values at the nodes, but also on the positive quadrature weights that we need to compute if they are not explicitly given.

This additional task is not necessary if we follow a different approach and, replacing the continuous scalar product (\ref{pr-cont}) by the following one
\begin{equation}\label{pr-dis}
<f,g>_N:=\sum_{i=1}^Nf(\xi_i)g(\xi_i),
\end{equation}
we consider the corresponding orthogonal polynomial projections, namely the least squares polynomials $\tilde\S_nf$ defined by
\[
\sum_{i=1}^N[f(\xi_i)-\tilde\S_nf(\xi_i)]^2=\min_{P\in\PP_n}
\sum_{i=1}^N[f(\xi_i)-P(\xi_i)]^2,
\]
and the following means of least squares polynomials
%\begin{equation}\label{meanLS}
\[
\tilde V_n^mf(\x)=\frac 1{2m+1}\sum_{r=n-m}^{n+m}\tilde\S_rf(\x),\qquad \qquad m=\lfloor \theta n\rfloor,\qquad 0<\theta<1.
\]
%\end{equation}
In this paper we are going to show that instead of $L_nf$ and $V_n^mf$ we can also use $\tilde\S_nf$ and $\tilde V_n^mf$ respectively, since they provide a comparable approximation w.r.t. the uniform norm.

More precisely, we point out that, unlike hyperinterpolation, least squares polynomial approximation does not require to know any quadrature weight, neither any quadrature rule is indeed necessary for its definition. Nevertheless, if we assume that (\ref{quad}) holds with $\mu=2n$, then, analogously to (\ref{LC-Fou}), we get (cf. Theorem \ref{th-ls})
\[
\|\tilde\S_n\|\sim\sqrt{n},
\]
which means that least squares similarly to hyperinterpolation polynomials  satisfy the following error estimate
\begin{equation}\label{hyper-err}
E_{n}(f)\le\|f-\tilde\S_nf\|_\infty\le  \C \sqrt{n} E_{n}(f),\qquad \C\ne\C(n,f).
\end{equation}
Moreover, as regards the means $\tilde V_n^mf$ of least squares polynomial approximants, we prove that for arbitrary $0<\theta<1$ and $m=\lfloor \theta n\rfloor$, similarly to (\ref{LC-VP}), we have (cf. Theorem \ref{th-VPls} )
\[
\sup_{n}\|\tilde V_n^m\|<\infty,
\]
so that $\tilde V_n^mf$ provides a near--best approximation of $f$ analogous to the filtered hyperinterpolation polynomials $V_n^mf$, satisfying the following error estimate
\begin{equation}\label{VP-err}
E_{n+m}(f) \le\|f-\tilde V_n^mf\|_\infty\le  \C E_{n-m}(f),\qquad m=\lfloor \theta n\rfloor,\qquad \C\ne\C(n,f).
\end{equation}
A similar result has been previously obtained in \cite{W3} only for the special case of spherical design configurations of points and for suitable sequences of filter coefficients defining delayed means of the following kind
$ \sum_{r=n-m}^{n+m}d_r\tilde\S_rf$, where $\sum_{r=n-m}^{n+m}d_r=1$.
The extension of (\ref{VP-err}) to these more general means will be given in a forthcoming paper by the authors.

The paper is organised as follows. In Section 2 we state our main results, illustrated by several numerical experiments.
The  proofs are given in Section 3 together with some auxiliary results, which can be of interest also in other contexts. Finally Section 4 is devoted to the conclusions.
\section{Main results}
Let $f$ be a given arbitrary function $f:\SS^2\rightarrow \RR$ such that $\|f\|_\infty<\infty$,  and suppose that we know its sampled values at the point set $X_N:=\{\xi_1,\ldots,\xi_N\}\subset\SS^2$.

As announced in the previous section, for our results we assume that the points in $X_N$ are nodes of a suitable positive weighted quadrature rule.

Let $|T|$ denote the cardinality of a set $T\subset\SS^2$, then the following sufficient conditions for the existence of positive weighted quadrature rules have been proven in \cite[Th. 3.1]{r943}.
\begin{theorem}\label{th-BD}
There exists a positive weighted quadrature rule of degree of exactness $4n$ based on the nodes $X_N$, if there exist suitable constants $a\ge 1$ and $\delta>0$ (independent of $n,N$) such that we have
\begin{equation}\label{hp}
\SS^2=\bigcup_{i=1}^N c\left(\xi_i, \frac \delta{n}\right),\qquad\mbox{and}\qquad\max_{1\le i\le N}\left|\left\{\xi_j\in X_N : \ d(\xi_j,\xi_i)\le \frac\delta{n}\right\}\right|\le a,
\end{equation}
where $c(\xi,\epsilon):=\{\x\in\SS^2 \::\: \ d(\xi,\x)\le\epsilon\}$, $d(\x,\y):=\arccos (\x\cdot \y)$ denotes the geodesic distance, and $\x\cdot\y$ is the Euclidean scalar product in $\RR^3$ of $\x,\y\in\SS^2$.
\end{theorem}
We point out that in \cite{r943} concrete bounds for the constants in (\ref{hp}) are not given, so that from the previous theorem we cannot derive any hint regarding the existence and the maximum degree of exactness of a positive quadrature rule based on the nodes in $X_N$. An attempt to find concrete bounds for the previous constants can be found in \cite{r938}, but the resulting theoretical bounds are too large \cite{W2} and the question remains an open problem to be further investigated.

Thus in the sequel we are going to assume that the points in $X_N$ are nodes of the positive quadrature rule (\ref{quad}) with $\mu=2n$ or $\mu=4n$. According to \cite[Th. 4.1]{W-Re1} this assumption implies the first condition in (\ref{hp}) is satisfied for some $\delta>0$ , while, at our knowledge, the second condition in (\ref{hp}) does not  necessarily follow from the existence of the quadrature rule.
Nevertheless if we assume that it holds, then the next result can be easily derived from \cite[eq. (3.6)]{r943}
\begin{theorem}\label{th-Marci}
Let $X_N=\{\xi_1,\ldots,\xi_N\}$  and $n\in\NN$ be such that
\begin{equation}\label{hp-Marci}
\max_{1\le i\le N}\left|\left\{\xi_j\in X_N : \ d(\xi_j,\xi_i)\le \frac\delta{n}\right\}\right|\le a
\end{equation}
holds with $\delta>0$  and $a\in\NN$ independent of $N$ and $n$. Then we have
\begin{equation}\label{Marci-1}
\frac 1{n^2}\sum_{i=1}^N|Q(\xi_i)|\le \C\int_{\SS^2}|Q(\x)|d\sigma(\x),\qquad \forall Q\in\PP_{n}, \qquad \C\ne\C(n,N,Q).
\end{equation}
\end{theorem}
The Marcinkiewicz type inequality (\ref{Marci-1}) is the second ingredient we need in stating our results.
%We remark that it may hold even if (\ref{hp-Marci}) is not satisfied (to check if this is for instance the case when $X_N$ are the nodes of Gauss--Legendre tensor product rules).
We remark  that if (\ref{Marci-1}) holds for all $Q\in\PP_n$, then it holds also for all $Q\in\PP_{ln}$ with $l\in\NN$ (in which case the constant $\C$ depends on $l$), namely (\ref{Marci-1}) implies that (see e.g. \cite[Th. 2.1]{W-Dai})
\begin{equation}\label{Marci}
\frac 1{n^2}\sum_{i=1}^N|Q(\xi_i)|\le \C\int_{\SS^2}|Q(\x)|d\sigma(\x),\qquad \forall Q\in\PP_{ln}, \quad l\in \NN, \qquad \C\ne\C(n,N,Q).
\end{equation}
To illustrate the theoretical results, we will take extremal systems of points as defined and computed in  \cite{W-pts} and call these point sets of type 1.
We recall that the point set $X_N=\{\xi_1,\ldots,\xi_N\}$ is said to be extremal if $N=\dim\PP_n$ and it maximizes the Vandermonde determinant
\[
\Delta(\xi_1,\ldots,\xi_N):=\det (\Phi_i(\xi_j))_{i,j=1}^N,
\]
where $\{\Phi_1,\ldots,\Phi_N\}$ is a basis of $\PP_n$.
It is known that extremal system of nodes are independent of the choice of the polynomial basis, and they support
%When the number of points $N$ is equal to $(n+1)^2$, then the point set supports
a positive quadrature rule of degree of exactness $n$.

\begin{figure}
\includegraphics[scale=0.5]{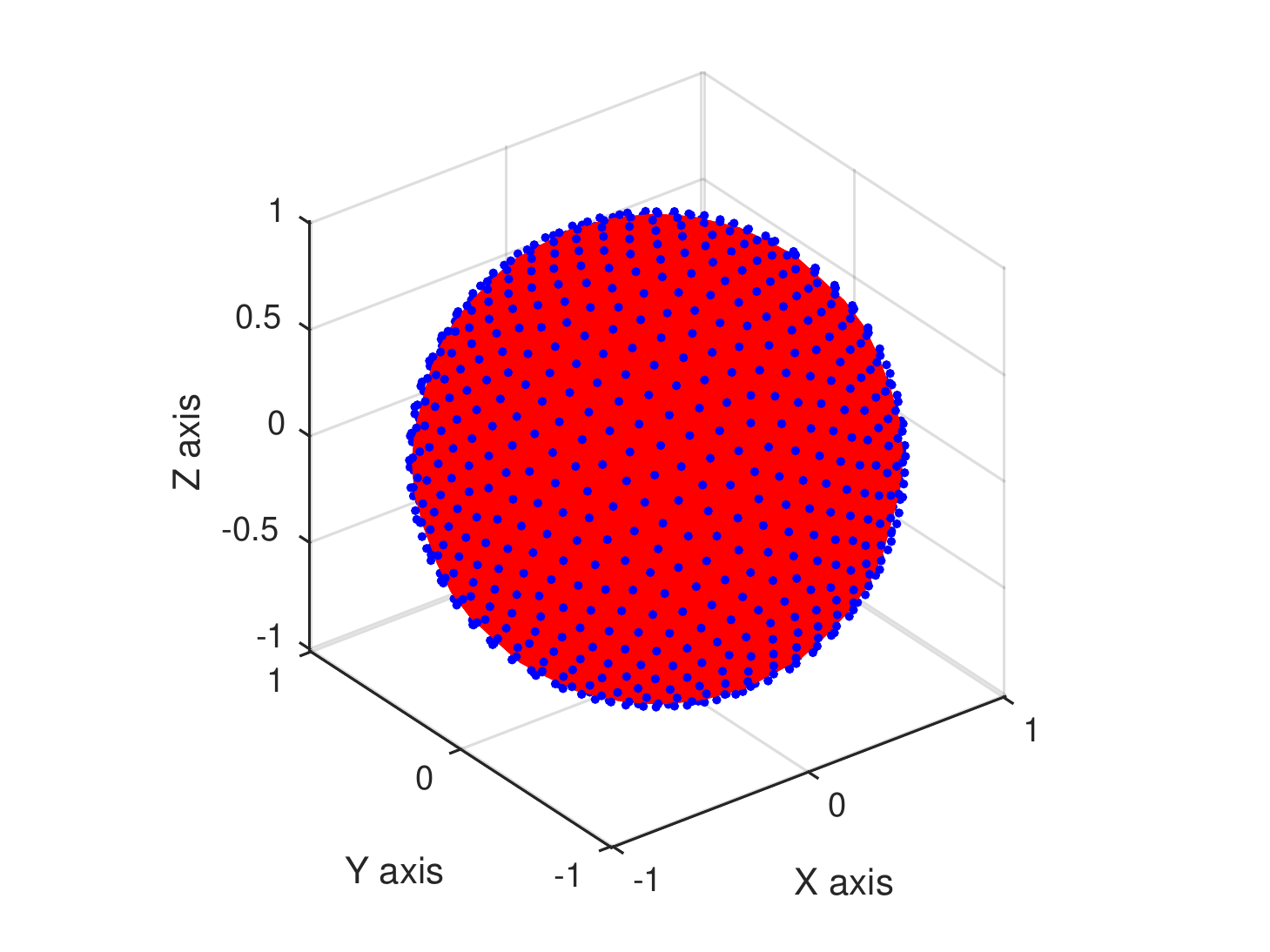}
\includegraphics[scale=0.5]{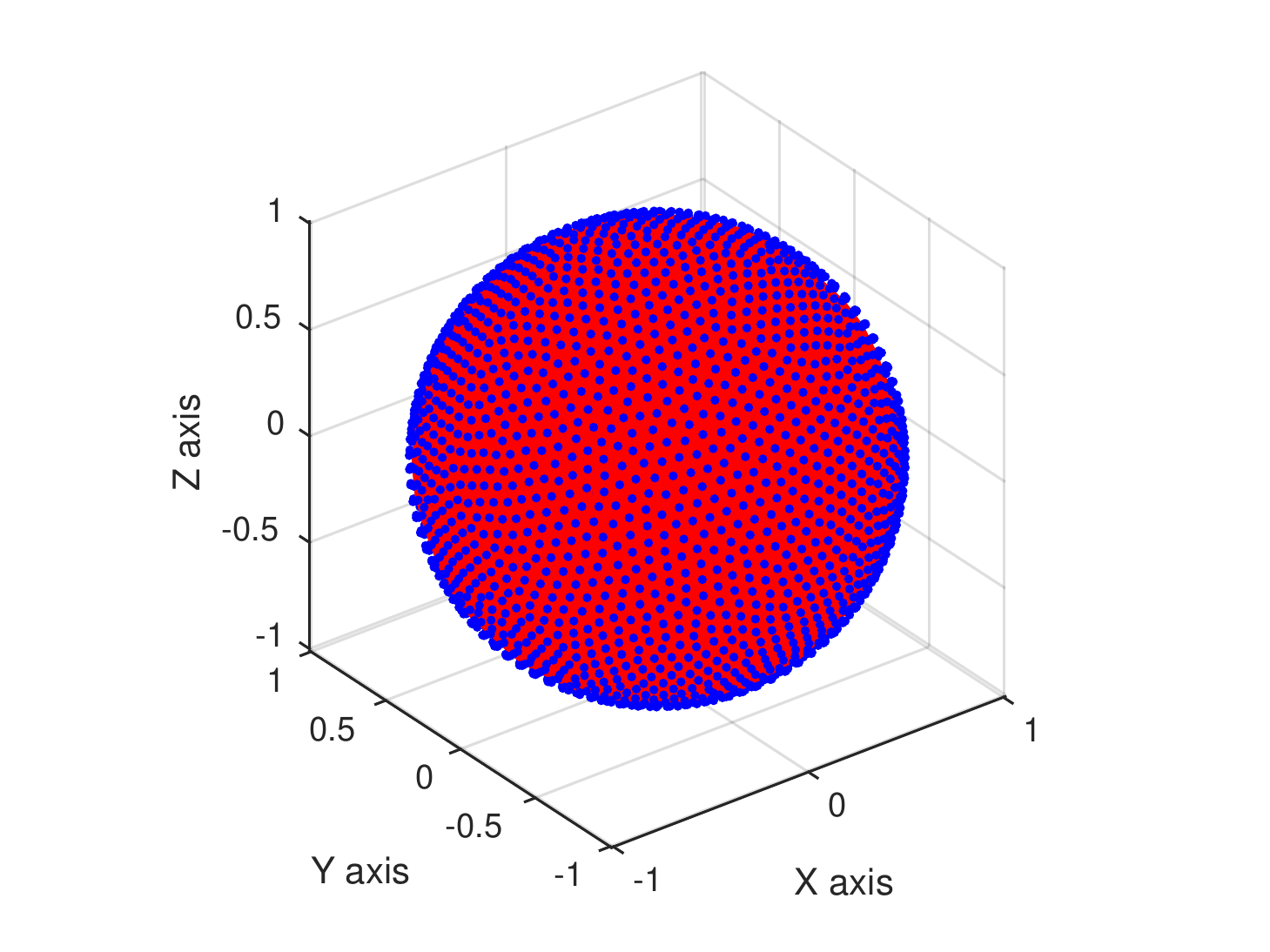}
\caption{Examples of the extremal systems of points related to degree of exactness $n=30$ (left) and $n=50$ (right).\label{fig001}}
\end{figure}

\begin{figure}
\includegraphics[scale=0.5]{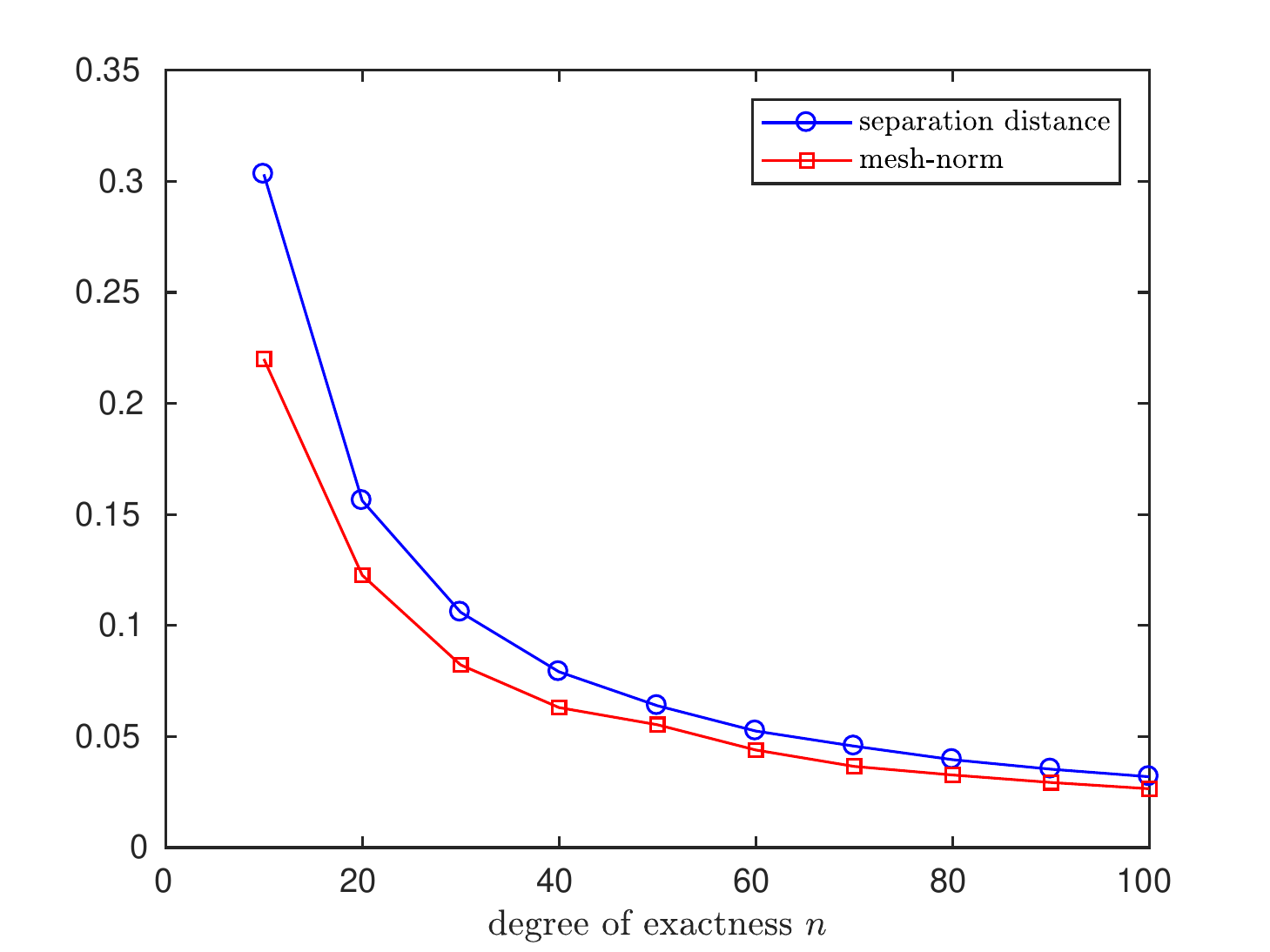}
%\includegraphics[scale=0.4]{fig012.pdf}
% \\
% \includegraphics[scale=0.5]{fig013.pdf}
%\includegraphics[scale=0.4]{fig014.pdf}
% \\
\includegraphics[scale=0.5]{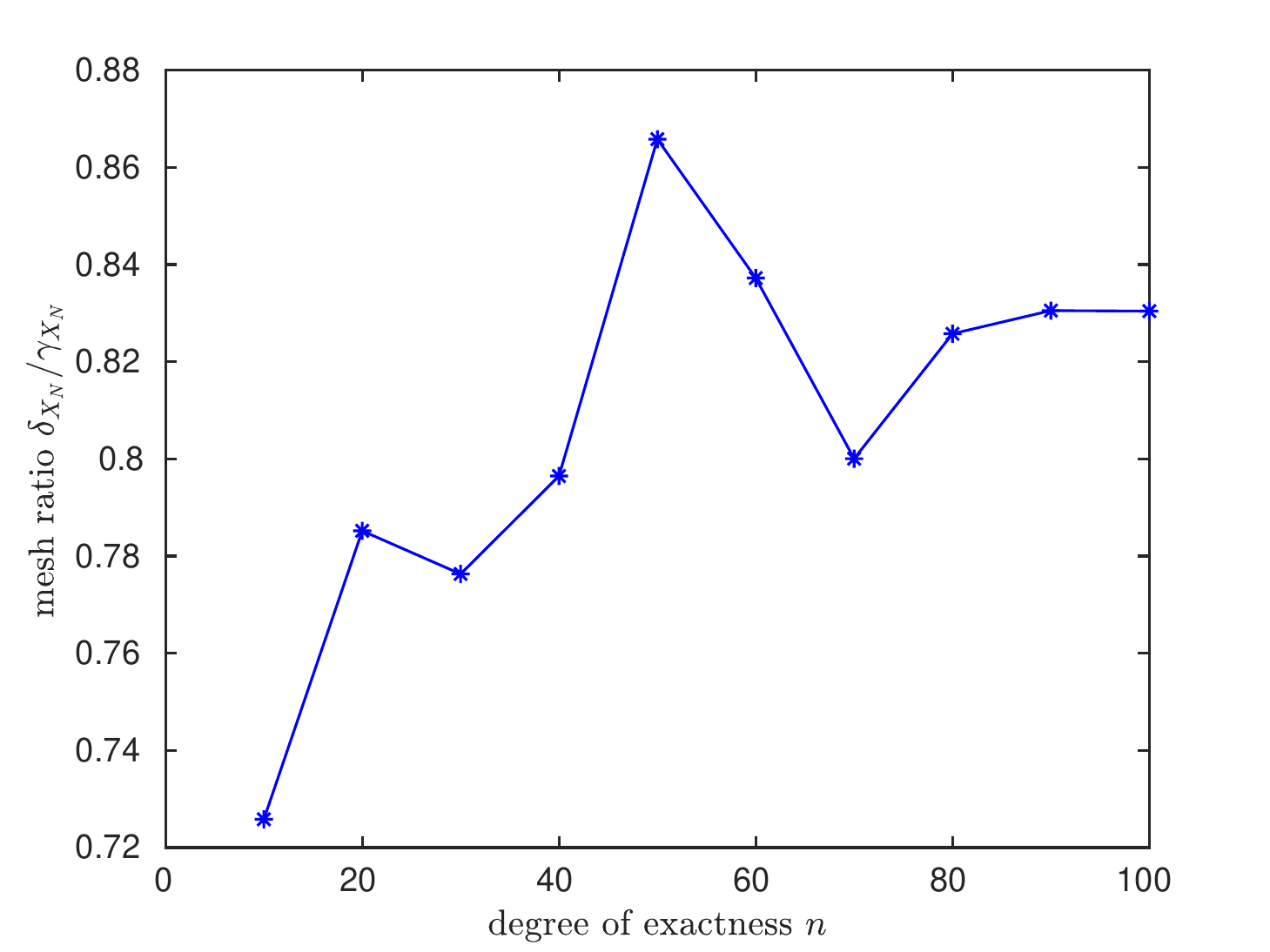}
\caption{The values of the separation distance $\gamma_{X_N}$ and the mesh-norm $\delta_{X_N}$ (left) and the mesh ratio $\delta_{X_N} / \gamma_{X_N}$ (right) for  extremal systems of $N=(n+1)^2$ points related to degrees of exactness $n=10,20,\ldots,100$.\label{fig011}}
\end{figure}

Let us now look at the different properties of point sets of type 1. To get an idea of how such a point set looks like,
we show two examples in Figure~\ref{fig001}, one for degree of exactness $n=30$ and one for degree $n=50$.
Note that the corresponding number of points is $N_{30} = 961$ and $N_{50} = 2601$.
In \cite[Section 5]{W-pts} the geometrical properties of the point sets of type 1 have been investigated by considering the mesh norm $\delta_{X_N}$ and separation distance $\gamma_{X_N}$ defined by (see e.g. (\cite{W-Po})
\begin{eqnarray}\label{delta}
\delta_{X_N}&:=&\max_{\x\in\SS^2}\min_{1\le i\le N} d(\x,\xi_i),\\
\label{gamma}
\gamma_{X_N}&:=&\min_{ i\ne j}d(\xi_i,\xi_j),
 \end{eqnarray}
It turns out \cite[Th. 5.1]{W-pts} that when $X_N$ is the point set of type 1, then we have
\[
\gamma_{X_N}\ge \frac \pi{2n},
\]
so that Theorem \ref{th-Marci} assures that (\ref{Marci}) holds.
Let us now consider each of the degrees of exactness $n=10,20,\ldots,100$ and compute the values of $\delta_{X_N}$ and $\gamma_{X_N}$ as
defined in (\ref{delta}) and (\ref{gamma}), respectively.
To estimate the mesh norm $\delta_{X_N}$, instead of taking the maximum over  the set of all
points of the sphere, the maximum is computed over a point set with a number of points considerably larger than the number of points for which we want
to approximate the mesh norm. To this end we consider the ``spiral points'' as defined in  \cite{HardMichSaff2016} which we call point sets of the second type.  These can be computed very efficiently and seem to be uniformly distributed
over the unit sphere where each of the points seems to be well separated from the others. They are also mentioned in \cite{r940}. These ``spiral points'' were generalized by Bauer \cite{Bauer2000} and called ``generalized spiral points'' in the overview paper  \cite{RakhSaffZhou1995}.
To estimate $\delta_{X_N}$, we take a point set of the second type having $16$ times more points.
The results are shown in Figure~\ref{fig011}.
%A good point distribution on the sphere is one where, for a given number of points,
%the distance between the points is as large as possible.
%In this respect, the first type of point sets does not score well as indicated in Figure~\ref{fig011}.

Now let us consider the least squares polynomial $\tilde\S_n f$ based  on the point set $X_N$, i.e., let
\begin{equation}\label{LS-min}
\sum_{i=1}^N[f(\xi_i)-\tilde\S_nf(\xi_i)]^2=\min_{P\in\PP_n}
\sum_{i=1}^N[f(\xi_i)-P(\xi_i)]^2.
\end{equation}
Taking into account that the map $\tilde\S_n:f\rightarrow \tilde\S_nf\in\PP_n$ is a projection onto $\PP_n$, i.e.,
\begin{equation}\label{inva-ls}
\tilde\S_nP=P,\qquad \forall P\in\PP_n,
\end{equation}
we easily get
\[
E_{n}(f)\le\|f-\tilde\S_nf\|_\infty\le  \left(1+\|\tilde\S_n\|\right) E_{n}(f),
\]
where $E_n(f)$ is the error of best uniform polynomial approximation (\ref{En}) and
\begin{equation}\label{LC-ls-def}
\|\tilde\S_n\|:=\sup_{g\ne 0}\frac{\|\tilde\S_ng\|_\infty}{\|g\|_\infty}.
\end{equation}
The constants in (\ref{LC-ls-def}) are usually called Lebesgue constants of $\tilde\S_n$ and their behaviour may strongly influence the quality of the approximation. With regard to this, it is known \cite{W-Da} that the projection onto $\PP_n$ having minimal Lebesgue constants is the Fourier orthogonal projection w.r.t. the scalar product (\ref{pr-cont}), which can be written in compact form as follows \cite{b509}
 \[
 \S_nf(\x)= \frac 1{2\pi}\int_{\SS^2}K_n(\x\cdot\y)f(\y)d\sigma(\y),\qquad \x\in\SS^2,
 \]
where $K_n(t)=\frac{(n+1)}2P_n^{(1,0)}(t)$, $\forall t\in [-1,1]$, and $P_n^{(1,0)}(t)$  is the Jacobi polynomial of degree $n$ associated with the weight $v^{(1,0)}(t)=(1-t)$ and normalized so that $P_n^{(1,0)}(1)=(n+1)$.

On the other hand, it is known that the Lebesgue constants of $\S_n$ grow with $n$ according to \cite{W-Da}
\begin{equation}\label{min-norm}
 \|\S_n\|\sim \sqrt{n}.
\end{equation}
By the next theorem, we state that also the Lebesgue constants of least squares projections have this minimal growth .
\begin{theorem}\label{th-ls}
Let $n\in\NN$ and let the set $X_N=\{\xi_1,\ldots,\xi_N\}\subset\SS^2$ be such that (\ref{Marci-1}) holds. If $X_N$ provides the nodes of a positive weighted quadrature rule of degree of exactness $2n$ (i.e., if (\ref{quad}) holds with $\mu=2n$ and $\lambda_i>0$),
then we have
\begin{equation}\label{Leb-LS}
\|\tilde\S_n\|\sim\sqrt{n}.
\end{equation}
\end{theorem}
In proving this theorem, a fundamental role is played by the existence of the positive weighted  quadrature rule (\ref{quad}) with $\mu=2n$. This condition also allows us to discretize the Fourier projection $\S_n$ obtaining the following  hyperinterpolation polynomial projection
\begin{equation}\label{hyper}
L_{n}f(\x):=\frac 1{2\pi}\sum_{i=1}^N\lambda_{i}f(\xi_i)K_{n}(\xi_i\cdot \x),\qquad \x\in\SS^2,
\end{equation}
which satisfies $\|L_n\|\sim \sqrt{n}$ too \cite{W4, W-Re}.
Nevertheless, we point out that unlike hyperinterpolation, the effective computation of the quadrature weights is not necessary for getting least squares polynomials.

\cmarc{
In order to illustrate the theoretical results, in the next experiment, we investigate the behaviour of the Lebesgue constant of both least squares and hyperinterpolation polynomials of degree $n$ related to the point set of type 1. To this end, we'll estimate the Lebesgue constant of the corresponding operator by taking a larger point set  of type 2 containing $4$ times the number of  points of type 1.
Figure~\ref{fig021} shows the results. The circles and squares indicate the Lebesgue constant for the discrete least squares operator and the hyperinterpolation operator, respectively, when we take for degree $n$ on the horizontal axis
the corresponding point set of type 1 related to the degree of exactness $2n$, i.e., having $N = (2n+1)^2$ points.
Moreover, the solid line and the dashed line indicate the Lebesgue constant for degree $n$ but where the point set is the point set of type 1
having $N = (2 \cdot 60 +1)^2$ points.
\begin{figure}
\includegraphics[scale=0.8]{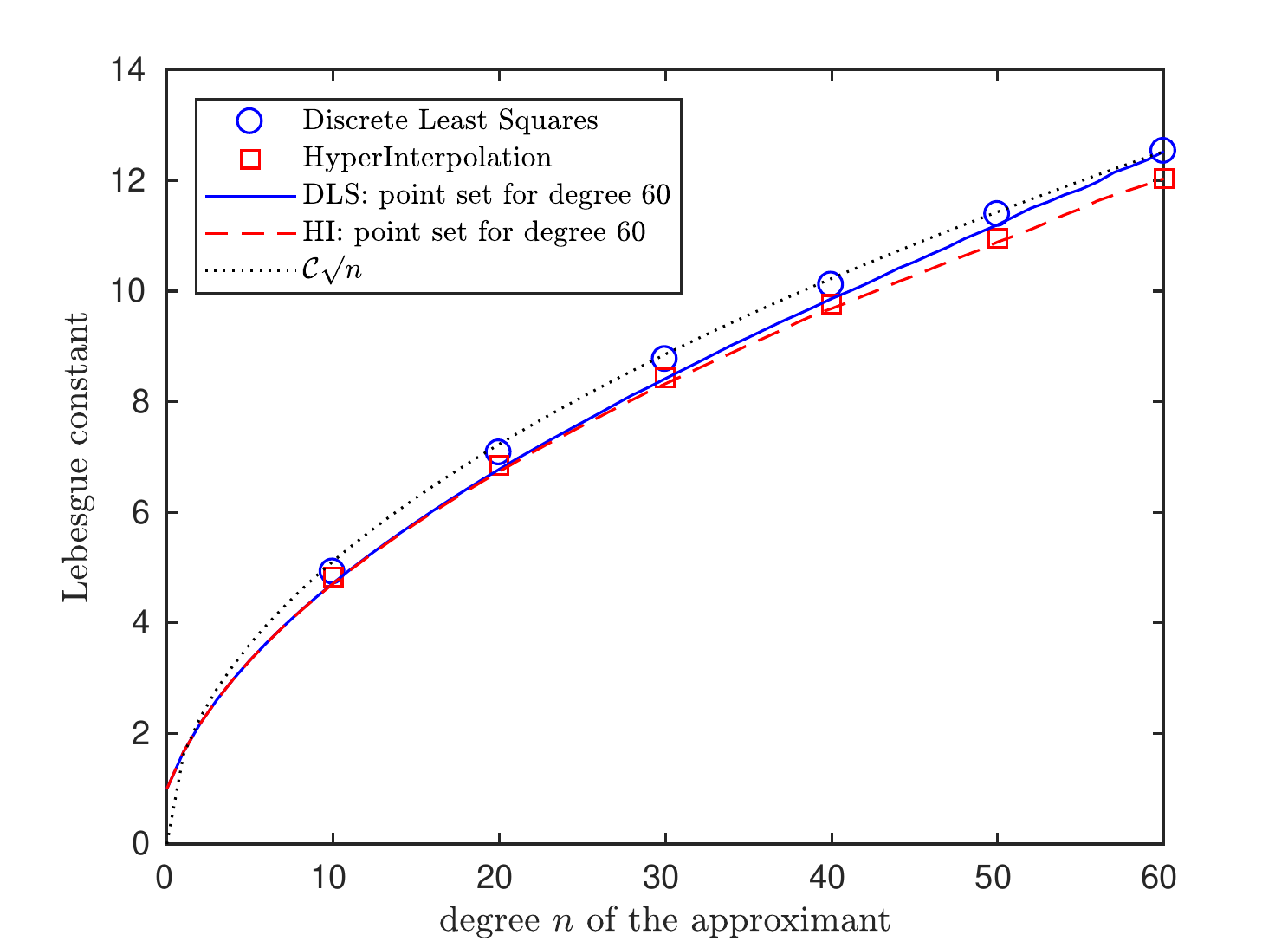}
%\\
%\includegraphics[scale=0.4]{fig023.pdf}
\caption{The values of the Lebesgue constants of the least squares operator $\tilde\S_n$ and the hyperinterpolation operator $L_n$ for the degrees $n=10,20,\ldots,60$ with corresponding point set of type 1 having $N = (2n+1)^2$ points, and for the degrees $n\le 60$
using the point set of type 1 having $N = (2\cdot 60+1)^2$ points.\label{fig021}}
\end{figure}
The picture shows that taking a larger point set only results in a slight decrease of the Lebesgue constant
while the Lebesgue constant for hyperinterpolation is slightly less than the one for discrete least squares approximation.
However, in computing the least squares approximant we do not need the weights of the quadrature rule.
%the point sets of the type considered here can not be used
%for polynomial interpolation, i.e., the corresponding generalized Vandermonde matrix is singular.
%However, once the degree of the least squares approximant is $10$ percent smaller than the degree
%of the point set, the $\infty$-norm is reasonably small as is shown in Figure~\ref{fig031}.
The optimal Lebesgue constant for polynomial projections on the unit sphere  grows as
 $\sqrt{n}$ (see, e.g., \cite[p.~40]{b509}).
This curve is also shown in Figure~\ref{fig021} (dotted line).
%\begin{figure}
%\includegraphics[scale=0.8]{fig031.pdf}
%\caption{The values of the $\infty$-norm of the least squares operator for the point sets related to degrees $n=10,20,\ldots,60$ where the degree of the least squares approximant is limited to $0.9$ times the degree of the point set.\label{fig031}}
%\end{figure}
%
}

In the case that $f$ is a continuous function (i.e., $f\in C(\SS^2)$) we recall that Weierstrass theorem extends to the sphere and it assures $f$ can be uniformly approximated by polynomials with the desired precision, namely
\[
\lim_{n\rightarrow\infty}E_n(f)=0,\qquad\forall f\in C(\SS^2),
\]
holds, the rate of convergence depending on the degree of smoothness of the function $f$ (see e.g. \cite{b509}).

On the other hand, the estimate (\ref{Leb-LS}) does not assure that, as $n\rightarrow\infty$, the sequence of least squares polynomials $\tilde\S_nf$ uniformly converges to $f$ for any continuous function $f$ on the sphere, but more smoothness properties are required on $f$ for getting the uniform convergence, according to the estimate
\begin{equation}\label{err-ls}
\|f-\tilde\S_nf\|_\infty\le  \C \sqrt{n} E_{n}(f),\qquad \C\ne\C(n,f).
\end{equation}
In order to get a sequence of discrete approximation polynomials which uniformly converges to any continuous function $f\in C(\SS^2)$, we are going to consider the following de la Vall\'ee Poussin type means of least squares polynomials
\begin{equation}\label{meanLS}
\tilde V_n^mf(\x)=\frac 1{2m+1}\sum_{r=n-m}^{n+m}\tilde\S_rf(\x),\qquad m=\lfloor \theta n\rfloor, \qquad \x\in\SS^2.
\end{equation}
Note that the map $\tilde V_n^m:f\rightarrow \tilde V_n^mf\in \PP_{n+m}$ is a polynomial quasi--projection, i.e., we have
\begin{equation}\label{inva-VPls}
\tilde V_n^mP=P,\qquad \forall P\in\PP_{n-m},
\end{equation}
which implies, by standard arguments, the following error estimates
\[
E_{n+m}(f)\le\|f-\tilde V_n^mf\|_\infty\le  \left(1+\|\tilde V_n^m\|\right) E_{n-m}(f),\qquad m=\lfloor \theta n\rfloor .
\]
The uniform boundedness of the Lebesgue constants $\|\tilde V_n^m\|$, is stated by the next theorem.
\begin{theorem}\label{th-VPls}
Let $n\in\NN$ and  $X_N=\{\xi_1,\ldots,\xi_N\}\subset\SS^2$ be such that (\ref{Marci-1}) holds. If the points in $X_N$ are nodes of a positive weighted quadrature rule of degree of exactness $4n$ (i.e., if (\ref{quad}) holds with $\mu=4n$ and $\lambda_i>0$), then for arbitrary $0<\theta<1$, and $m=\lfloor \theta n\rfloor$, we have
\begin{equation}\label{LC-VPls}
\sup_{n}\|\tilde V_n^m\|<\infty.
\end{equation}
\end{theorem}
We remark that if the points set $X_N$ supports a positive weighted quadrature rules of degree $4n$, then we can apply this formula in order to get the discrete counterpart of the mean in (\ref{VP-cont}), namely
\begin{equation}\label{filtered}
V_n^mf(\x):=\sum_{i=1}^N\lambda_i f(\xi_i)\left[\frac 1{2m+1}\sum_{r=n-m}^{n+m}K_r(\x\cdot\xi_i)\right],\qquad m=\lfloor\theta n\rfloor,\qquad 0<\theta<1.
\end{equation}
This polynomial falls in the class of filtered hyperinterpolation polynomials firstly introduced in \cite{r951}. It also satisfies (\ref{inva-VPls}) and has uniformly bounded Lebesgue constants \cite{W1}, but using $\tilde V_n^mf$ instead of $V_n^mf$ we do not need to compute the quadrature weights $\lambda_i$.

\cmarc{
In Figure~\ref{fig051}, the Lebesgue constants are plotted for $\tilde V_n^m$ (circles) and $V_n^m$ (squares).
For each value of $n = 5,10,15,\ldots,40$ the corresponding point set of type 1 related to degree of exactness $4n$ is taken, i.e.,
having $N = (4n+1)^2$ points. To get the estimates point sets of type 2 are considered having $4N$ points.
%the discrete de la Vall\'ee Poussin (VP) operator.
%Given the degree of the point set $\delta$, in the VP mean the highest degree of the approximant will be $n+m$.
%As we have seen before, when we would take $n+m = \alpha \delta$ with $\alpha \approx 1$ the $\infty$-norm of the least squares approximant
%of degree $n+m$ could be large. Hence, we have to choose $\alpha$ smaller than $1$.
The coupling between the parameters $n$ and $m$ is $m=\lfloor \theta n\rfloor$
with $0\le\theta \leq 1$. In the left subfigure we show the different lines for $\theta = 0.0, 0.1, \ldots, 1.0$ from top to bottom
while for the right subfigure we focus on $\theta = 0.1, 0.2, \ldots, 1.0$. Note that there is only a small difference between the
Lebesgue constants for the mean of least squares approximations and the filtered hyperinterpolation case.
\begin{figure}
\includegraphics[scale=0.5]{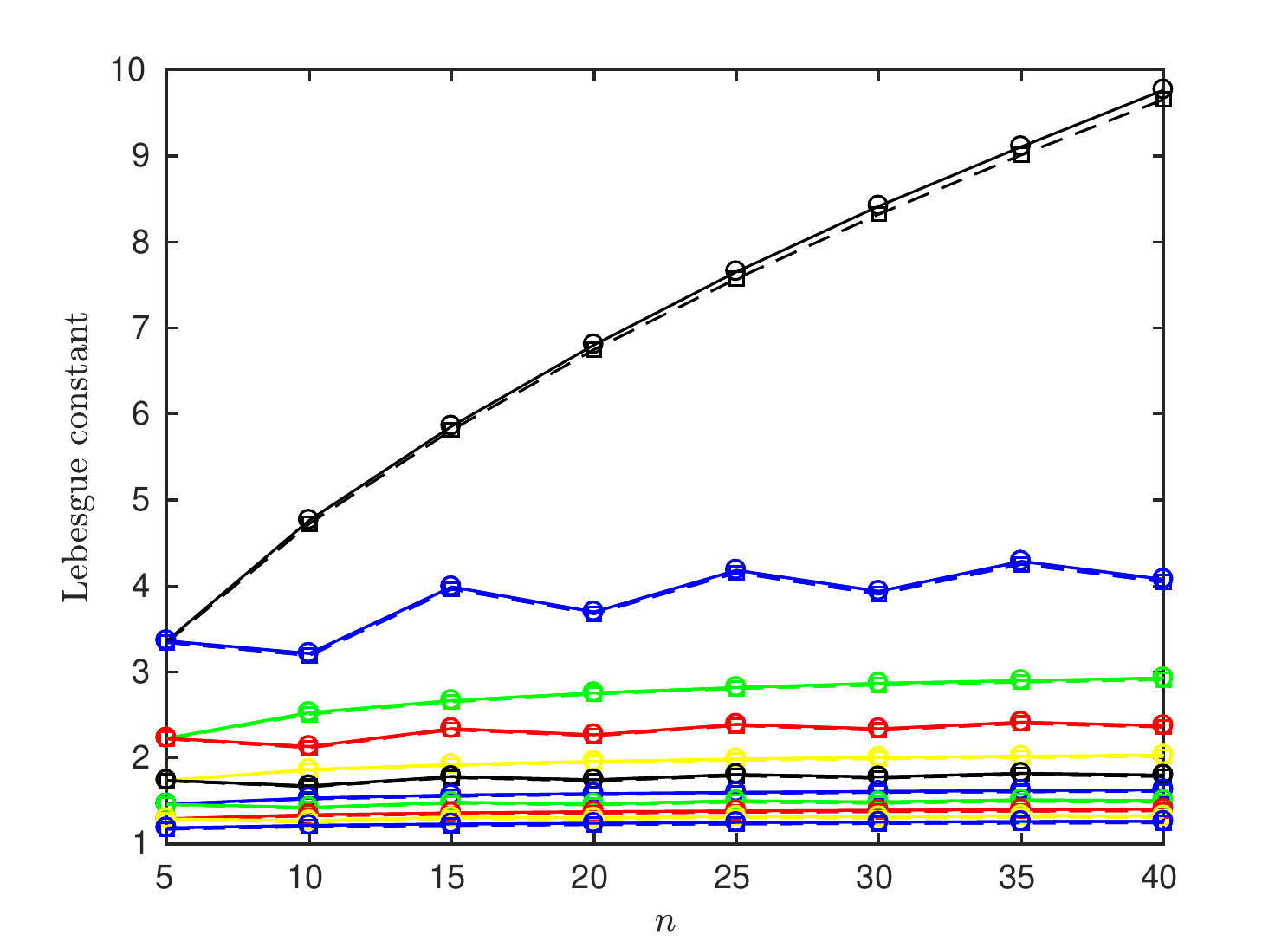}
\includegraphics[scale=0.5]{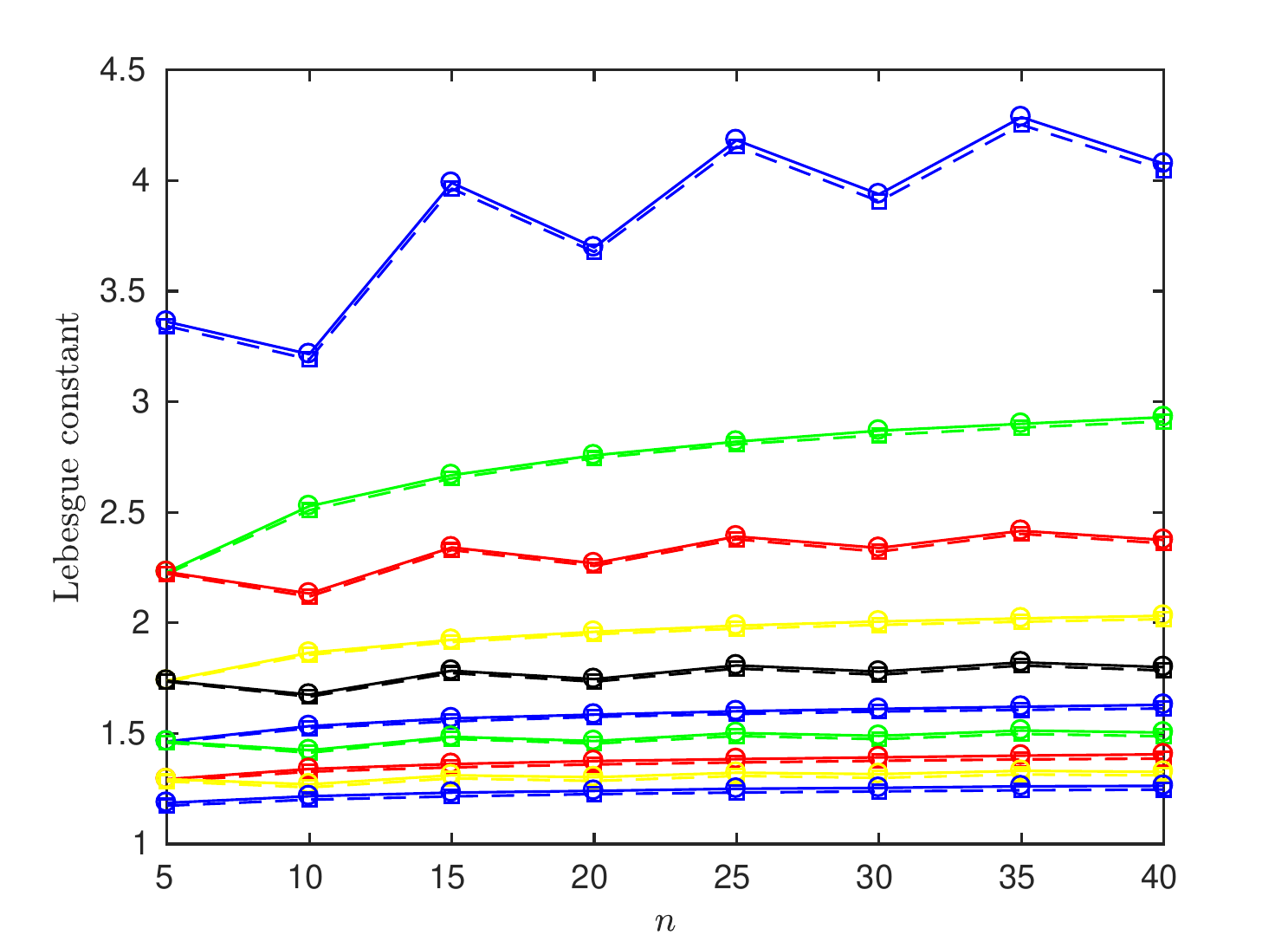}
\caption{The values of the Lebesgue constant of $\tilde V_n^m$ (circles) and $V_n^m$ (squares) in function of the degree $n$. At the left side,
the curves show the results for $\theta = 0.0,0.1,\ldots,1.0$ (from top to bottom)  while it starts from $\theta = 0.1$ at the right side.\label{fig051}}
\end{figure}

In the case of functions almost everywhere smooth, apart from some isolated points of singularity, it is known that the Gibbs phenomenon occurs by using least squares as well as hyperinterpolation polynomials. In the following experiment we show that similarly to filtered hyperinterpolation, this phenomenon can be strongly reduced if we consider the mean $\tilde V_n^mf$. Keeping $n$  unchanged (which is strictly related to the number $N$ of data points) we appropriately modulate the range of action $m$ of the mean by suitably varying the parameter $\theta$ from the limiting value $\theta=0$, which corresponds to simple least squares approximation (suggested for very smooth functions) to the limiting value $\theta =1$, which in practice corresponds to the Fej\'er mean of least squares polynomials.
We consider the function $f_2$ as defined in \cite{r951}.
For completeness we repeat the definition of this function.
\[
f_{2}(\x) = f_{1}(\x) + f_{cone}(\x), \qquad \x\in\SS^2,
\]
where
\begin{eqnarray*}
f_{1}(\x) &=& 0.75 \exp(-(9x - 2)^2/4 - (9y - 2)^2/4 - (9z - 2)^2/4) \\
 && +0.75 \exp(-(9x + 1)^2/49 - (9y + 1)/10 - (9z + 1)/10) \\
 && +0.5 \exp(-(9x - 7)^2/4 - (9y - 3)^2/4 - (9z - 5)^2/4) \\
 && -0.2 \exp(-(9x-4)^2-(9y-7)^2-(9z-5)^2),
\end{eqnarray*}
and
$$
f_{cone}(\x) = \left\{
\begin{array}{ll}
2 \left( 1 - \frac{d(\x_c,\x)}{r} \right) & \mbox{ if } d(\x_c,\x) \leq r, \\
0 & \mbox{ if } d(\x_c,\x) > r,
\end{array}
\right.
$$
with $r = \frac{1}{2}$ and $\x_c = (\frac{1}{2}, \frac{1}{2}, \frac{1}{\sqrt{2}})^T$.

In Figure \ref{fig080}, the function $f_{2}$ is shown as well as the error for the  $\tilde V_n^mf_2$ approximants with
$n=20$, $m = \lfloor \theta n \rfloor$ with $\theta = 0.0, 0.1$ and $0.2$.
The point set is of type 1 and supports a positive quadrature rule of degree of precision $4n$, i.e., it has $N = (4n+1)^2$ points.
The figure can be compared to Figure $4$ in~\cite{r951}.
\begin{figure}
\includegraphics[scale=0.25]{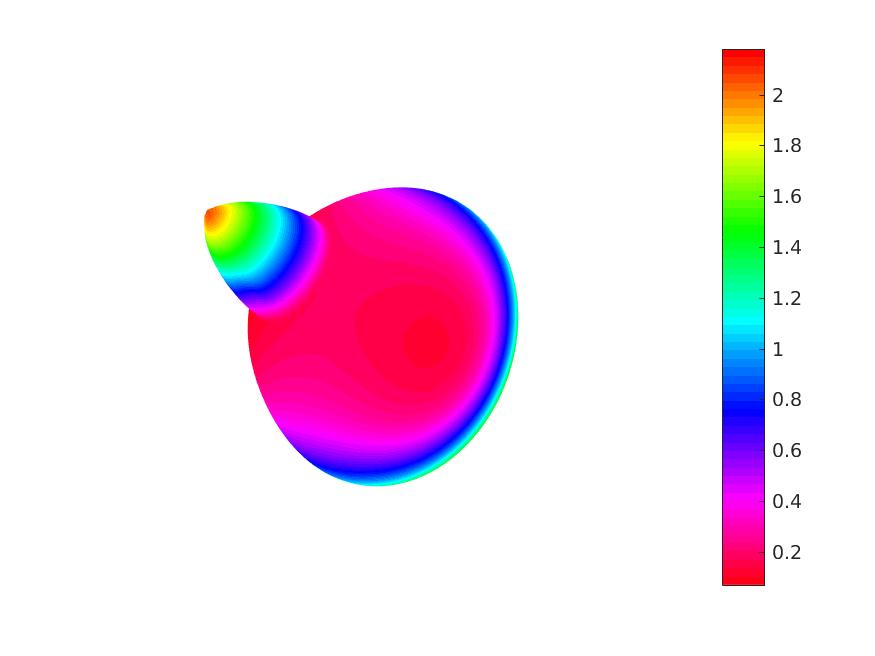}
\includegraphics[scale=0.25]{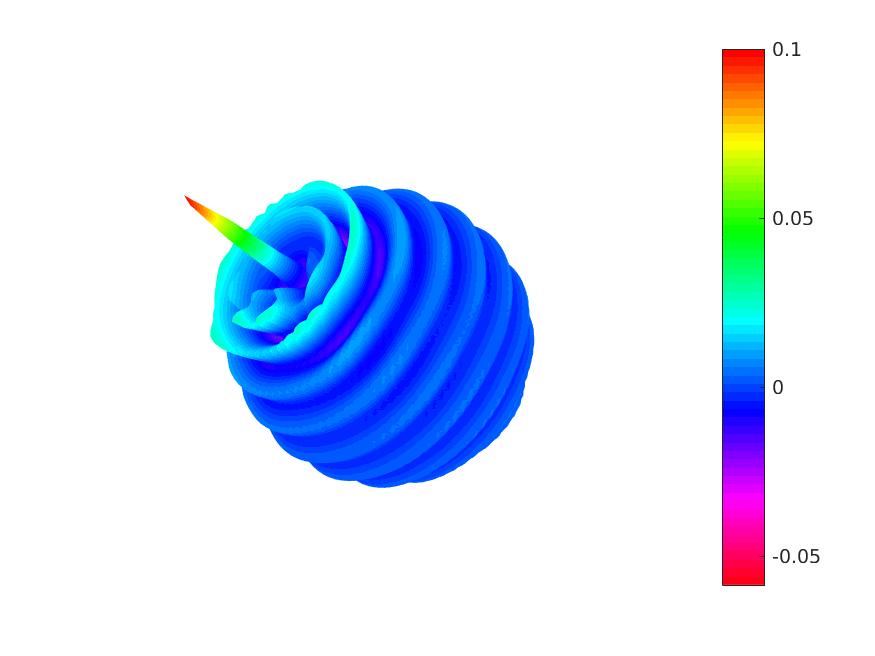} \\
\includegraphics[scale=0.25]{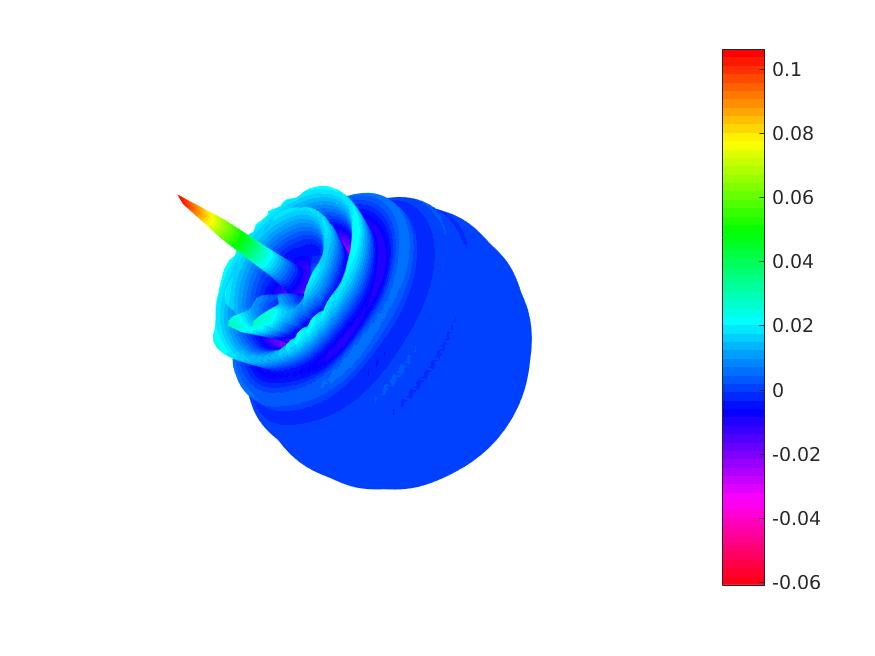}
\includegraphics[scale=0.25]{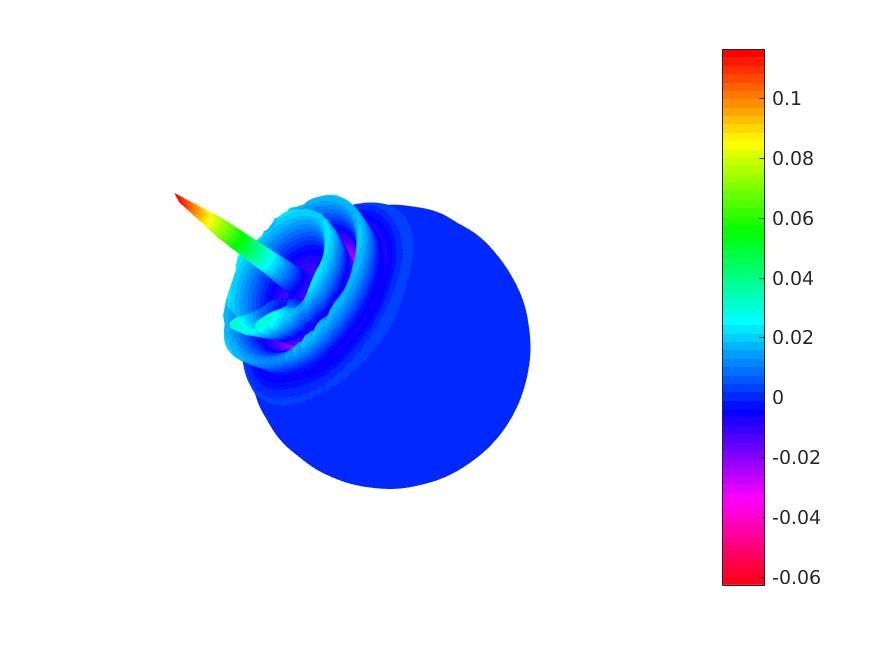}
\caption{In the top-left figure the function $f_{2}$ is shown. The plots in the top-right, bottom-left and bottom-right represent
the error $\|f_{2} - \tilde{V}_m^nf_{2} \|_\infty$ for $n = 20$, $m=\lfloor \theta n \rfloor$ and  $\theta = 0.0, 0.1, 0.2$.\label{fig080}}
\end{figure}
}

\section{Proofs and auxiliary results}
In order to prove the theorems of the previous section, we are going to state some preliminary results.

Firstly we recall that the space $\PP_n$ of spherical polynomials of degree at most $n$ has dimension $(n+1)^2$ and it is generated by spherical harmonics, which constitute an orthonormal basis  w.r.t. the scalar product (\ref{pr-cont}).

We refer to the literature (see, e.g., \cite{b509}) for  background information on spherical harmonics. Here we only recall that they are related to Legendre polynomials by an addition formula (cf. \cite[(1.6.7)]{b509}) so that Fourier orthogonal projections $\S_nf$ of any function $f$ into $\PP_n$ can be written as follows
\begin{equation}\label{Fourier}
\S_nf(\x)= \frac 1{2\pi}\int_{\SS^2}K_n(\x\cdot\y)f(\y)d\sigma(\y),\qquad \x\in\SS^2,
\end{equation}
with the reproducing kernel given by
\begin{equation}\label{Darboux}
K_n(t):=K_n(t,1)=\frac{(n+1)}2P_n^{(1,0)}(t),\qquad t\in [-1,1],
\end{equation}
where $K_n(t,s)$ denotes the Legendre--Darboux kernel of order $n$, and $P_n^{(1,0)}(t)$  is the Jacobi polynomial of degree $n$ associated with the weight $v^{(1,0)}(t)=(1-t)$ and normalized so that $P_n^{(1,0)}(1)=(n+1)$.

Obviously we have
\begin{equation}\label{inva-Fou}
P(\x)=\int_{\SS^2}P(\y)K_n(\x\cdot\y)d\sigma(\y),\qquad \forall P\in\PP_n,\qquad \forall\x\in\SS^2,
\end{equation}
and from (\ref{LC-Fou}), we deduce that
\begin{equation}\label{Leb-Fou}
\sup_{\x\in\SS^2}\left[\frac1{2\pi}\int_{\SS^2}|K_n(\x\cdot\y)|d\sigma(\y)\right]= \|\S_n\|\sim\sqrt{n}.
\end{equation}
Moreover, de la Vall\'ee Poussin mean (\ref{VP-cont}) can be written as follows
\begin{equation}\label{VP-cont-int}
\V_n^mf(\x)%=\frac 1{2m+1}\sum_{r=n-m}^{n+m}\S_rf(\x)
= \frac 1{2\pi}\int_{\SS^2}\left[\frac1{2m+1}
\sum_{r=n-m}^{n+m}K_r(\x\cdot\y)\right]f(\y)d\sigma(\y),\qquad \x\in\SS^2,
\end{equation}
and for all $m=\lfloor \theta n\rfloor$ with $0<\theta<1$, by (\ref{LC-VP}) we get
\begin{equation}\label{LC-VPclas}
\sup_{\x\in\SS^2}\left(\frac1{2\pi}
\int_{\SS^2}\left|\frac1{2m+1}\sum_{r=n-m}^{n+m}K_r(\x\cdot\y)\right|d\sigma(\y)\right)= \|\V_n^m\|\le  \C\ne\C(n).
\end{equation}
The next lemma generalizes this result.
\begin{lemma}\label{lem-gen}
Let $0<\theta<1$ and $m=\lfloor \theta n\rfloor$, with $n\in\NN$. Then for all $\x\in\SS^2$, we have
\begin{equation}\label{eq-lem-gen}
\sup_{n-m\le s\le n+m}\left(\int_{\SS^2}\left|\frac1{2m+1}\sum_{r=n-m}^{s}K_r(\x\cdot\y)\right|d\sigma(\y)\right)\le\C,\qquad \C\ne \C(n,\x).
\end{equation}
\end{lemma}
{\it Proof of Lemma \ref{lem-gen}}

Since it is easy to check that
\begin{equation}\label{int}
\int_{\SS^2}\varphi(\x\cdot\y)d\sigma(\y)=2\pi\int_{-1}^1\varphi(t)dt,
\qquad \forall\varphi\in L^1[-1,1],\quad \forall \x\in\SS^2,
\end{equation}
then, for all $\x\in\SS^2$ and for any $s=(n-m),\ldots,(n+m)$, we have
\[
\int_{\SS^2}\left|\frac1{2m+1}\sum_{r=n-m}^{s}K_r(\x\cdot\y)\right|d\sigma(\y)=2\pi
\int_{-1}^1\left|\frac1{2m+1}\sum_{r=n-m}^{s}K_r(t)\right|dt=: 2\pi I_s.%=2\pi\mathbb{V}_sg(1),
\]
In the case $s=n-m$, recalling that (see, e.g., \cite{r938, b210})
\begin{equation}\label{sup-Darboux}
\sup_{|t|\le 1}|K_n(t)|=K_n(1)=\frac{(n+1)^2}2,
\end{equation}
we get
\[
I_s=\int_{-1}^1\left|\frac1{2m+1}K_{n-m}(t)\right|dt\le \C\ne \C(n).
\]
For the case $n-m<s\le n+m$ we set for all $t\in [-1,1]$
\[
\mathbb{V}_sg(t)=\int_{-1}^1 g(u) \left(\frac1{2m+1}\sum_{r=n-m}^{s}K_r(t,u)\right)du,
\qquad\quad g(t):=\sign \left[\sum_{r=n-m}^{s}K_r(t)\right],
\]
and observe that $I_s=\mathbb{V}_sg(1)$.\newline
Then, taking into account that $m\sim n\sim s\sim (n-m)$ follows from our assumptions, we can apply \cite[Theorem 3.1]{Marc}, which yields
\[
|I_s|=|\mathbb{V}_sg(1)|\le \sup_{|t|\le 1}|\mathbb{V}_sg(t)|\le \C \sup_{|t|\le 1}|g(t)| = \C\ne \C(s,n,g),
\]
and the statement follows.
\Proofend

%\cmarc{The theorem is true also for $s = n-m$. This is because the factor $1/(2m+1)$ is still in the formula Otherwise, this would not be true?}

In stating our main results we always assume the existence of a positive quadrature rule (\ref{quad}).
By the next lemma we list some necessary conditions, which will be useful in the sequel.
\begin{lemma}\label{lem-quad}
If the quadrature rule (\ref{quad}) holds with positive weights $\lambda_i>0$, $i=1,\ldots,N$, and degree of exactness $\mu=2n$, then we have
\begin{equation}\label{li}
(n+1)^2<N,\qquad \mbox{and}\qquad \lambda_i\le\frac{4\pi}{(n+1)^2},\qquad i=1,\ldots,N.
\end{equation}
Moreover, for all $0<\rho<1$, we have
\begin{equation}\label{Marci-inf}
\|P\|_\infty\le\C\max_{1\le i\le N}|P(\xi_i)|,\qquad\forall P\in\PP_{\lfloor\rho n\rfloor},\qquad \C\ne\C(P,n,N).
\end{equation}
\end{lemma}
{\it Proof of Lemma \ref{lem-quad}}

In \cite{W-Re}[p. 274] it has been already observed that a necessary condition for the existence of a positive quadrature rule of degree of exactness $2n$, based on $N$ points, is $(n+1)^2<N$, which is the first bound in (\ref{li}).

In order to state the second bound in (\ref{li}),  for any $i=1,\ldots, N$, we apply (\ref{quad}) to the non negative polynomial $P(\x)=\left[\frac{K_n(\x\cdot\xi_i)}{K_n(1)}\right]^2\in\PP_{2n}$, and using (\ref{int}) and well--known properties of Legendre polynomials, we get
\begin{eqnarray*}
\lambda_{i}&=&
\lambda_{i}\left[\frac{K_n(\xi_i\cdot\xi_i)}{K_n(1)}\right]^2\le
\sum_{s=1}^N\lambda_{s}\left[\frac{K_n(\xi_s\cdot\xi_i)}
{K_n(1)}\right]^2=\int_{\SS^2}\left[\frac{K_n(\x\cdot\xi_i)}{K_n(1)}\right]^2d\sigma(\x)\\
&=& \frac{2\pi}{[K_n(1)]^2}\int_{-1}^1K_n^2(x)dx=\frac {2\pi}{K_n(1)}=\frac{4\pi}{(n+1)^2}.
\end{eqnarray*}
Finally, in order to prove (\ref{Marci-inf}), for any $n\in\NN$ and $0<\rho<1$ we set
\[
\theta:=1-\rho,\qquad m:=\lfloor\theta n\rfloor,
\]
so that we have $\lfloor\rho n\rfloor=n-\lfloor\theta n\rfloor=n-m$.

Then, for all $P\in\PP_{\lfloor\rho n\rfloor}=\PP_{n-m}$ and any $\x\in\SS^2$, by means of (\ref{inva-Fou}), (\ref{quad}), (\ref{li}), (\ref{Marci})  and (\ref{LC-VPclas}) we get
\begin{eqnarray*}
|P(\x)|&=& %|\V_n^mP(\x)|=
\left|\int_{\SS^2}P(\y)\left(\frac 1{2m+1}\sum_{r=n-m}^{n+m}K_r(\x\cdot\y)\right)d\sigma(\y)\right|\\
&=& \left|\sum_{i=1}^N\lambda_i P(\xi_i)\left(\frac 1{2m+1}\sum_{r=n-m}^{n+m}K_r(\x\cdot\xi_i)\right)\right|\\
&\le&\max_{1\le i\le N} |P(\xi_i)|\left(\sum_{i=1}^N\lambda_i\left|\frac 1{2m+1}\sum_{r=n-m}^{n+m}K_r(\x\cdot\xi_i)\right|\right)\\
&\le&\C \max_{1\le i\le N} |P(\xi_i)|\left(\frac 1{n^2}\sum_{i=1}^N\left|\frac 1{2m+1}\sum_{r=n-m}^{n+m}K_r(\x\cdot\xi_i)\right|\right)\\
&\le&\C \max_{1\le i\le N} |P(\xi_i)|\left(\int_{\SS^2}\left|\frac 1{2m+1}\sum_{r=n-m}^{n+m}K_r(\x\cdot\y)\right|d\sigma(\y)\right)\\
&\le&\C \max_{1\le i\le N} |P(\xi_i)| .
\end{eqnarray*}
\Proofend

For proving the theorems of the previous section, it is useful to derive an explicit expression for least squares polynomial approximants $\tilde \S_nf$.
To this aim, we recall that for all $n\in\NN$, by applying Gram--Schmidt orthogonalization process to the spherical harmonics basis of $\PP_n$, we can obtain a basis of $\PP_n$ orthonormal w.r.t. the discrete scalar product (\ref{pr-dis}). Hence, for all $n\in\NN$, we assume
$\PP_n= span \{I_r : r=1,\ldots,(n+1)^2\} $ with
\begin{equation}\label{ort}
<I_k, I_h>_N=\sum_{j=1}^NI_k(\xi_j)I_h(\xi_j)=\delta_{h,k}:=\left\{
\begin{array}{ll}
1 & h=k\\
0 & h\ne k
\end{array}\right.
,\qquad h,k=1,2,\ldots
\end{equation}
Defining $H_n(\x,\y)$ as
\begin{equation}\label{LS-ker}
 H_n(\x,\y):=\sum_{r=1}^{(n+1)^2}I_r(\x)I_r(\y),\qquad \x,\y\in\SS^2,
\end{equation}
it is easy to check that the least squares polynomial $\tilde\S_nf$ can be written in explicit form as follows
\begin{equation}\label{LS-sum}
\tilde\S_{n}f(\x)=\sum_{i=1}^Nf(\xi_i)H_{n}(\x, \xi_i),\qquad \x\in\SS^2.
\end{equation}
Consequently, by (\ref{inva-ls}) we get
\begin{equation}\label{inva-LSker}
P(\x)=\sum_{i=1}^NP(\xi_i)H_{n}(\x, \xi_i),\qquad \x\in\SS^2,\qquad \forall P\in\PP_n,
\end{equation}
and the Lebesgue constants are explicitly given by
\begin{equation}\label{LC-ls}
\|\tilde\S_n\|=\sup_{\x\in\SS^2} \left[\sum_{k=1}^N|H_{n}(\x, \xi_k)|\right].
\end{equation}
Finally let us prove that
\begin{equation}\label{ls-ker1}
|H_n(\xi_i,\xi_j)|\le\sum_{r=1}^{(n+1)^2}|I_r(\xi_i)I_r(\xi_j)|\le 1,\qquad \forall \xi_i,\xi_j\in X_N\qquad\forall N\ge (n+1)^2.
\end{equation}
Indeed, we note that the (rectangular) matrix $I := [I_k(\xi_h)]_{h=1,...,N}^{k=1,2,\ldots,(n+1)^2}$ formed by the orthonormal columns $[I_k(\xi_1),\ldots I_k(\xi_N)]^T$, $k=1,\ldots, (n+1)^2$, can be extended by additional columns to form a square orthogonal matrix $Q=[Q_{h,k}]_{h,k=1,...,N}$ with $Q_{h,k}=I_k(\xi_h)$, $k=1,2,\ldots,(n+1)^2$. Hence, we have
\begin{equation}\label{Q}
\sum_{k=1}^{{(n+1)^2}}|I_k(\xi_i)|^2\le \sum_{k=1}^{N}|Q_{i,k}|^2 =1,\qquad i=1,\ldots,N,\qquad \forall N\ge (n+1)^2,
\end{equation}
and consequently for all $n\in\NN$ s.t. $(n+1)^2\le N$ and for any pair of nodes $\xi_i,\xi_j\in X_N$, we get
\[
|H_n(\xi_i,\xi_j)|\le\sum_{r=1}^{(n+1)^2}|I_r(\xi_i)I_r(\xi_j)|\le \left(\sum_{r=1}^{{(n+1)^2}}|I_r(\xi_i)|^2\right)^\frac 12
\left(\sum_{r=1}^{{(n+1)^2}}|I_r(\xi_j)|^2\right)^\frac 12\le 1.
\]
Using the previous results, in the next subsections we are going to prove the theorems of Section 2.
\subsection{Proof of Theorem \ref{th-ls}}
By (\ref{LC-ls}) it is sufficient to prove that $\forall \xi\in\SS^2$ we have
\begin{equation}\label{tesi}
\sum_{k=1}^N|H_{n}(\xi, \xi_k)|\le\C\sqrt{n},\qquad \C\ne\C(n,N,\xi).
\end{equation}
We first state (\ref{tesi}) in the case that $\xi\in X_N$ and then we are going to extend this property to all $\xi\in\SS^2$.

Let $\xi\in X_N$ be arbitrarily fixed. Using (\ref{inva-Fou}) with $P(\x)=H_n(\xi,\x)$, and recalling that (\ref{quad}) holds with $\mu=2n$, we get
\[
H_n(\xi,\xi_k)=\int_{\SS^2}H_n(\xi,\y)K_n(\xi_k\cdot\y)d\sigma(\y)=\sum_{i=1}^N\lambda_i  H_n(\xi,\xi_i)
K_n(\xi_k\cdot\xi_i),\qquad k=1,\ldots,N.
\]
On the other hand,
supposing (without loss of generality) that the quadrature weights are such that
\begin{equation}\label{l}
\lambda_0:=0<\lambda_1\le \lambda_2\le\ldots\le\lambda_N,
\end{equation}
if we apply the following summation by parts formula (with $a_i=\lambda_i$)
\begin{equation}\label{sum-part}
\sum_{i=1}^Na_ib_i=a_1\sum_{i=1}^N b_i+\sum_{i=2}^{N}(a_{i}-a_{i-1})\sum_{j=i}^Nb_j,
\end{equation}
then, by (\ref{li}), we get
\begin{eqnarray*}
\sum_{k=1}^N|H_n(\xi,\xi_k)|&=&\sum_{k=1}^N\left|\sum_{i=1}^N\lambda_i  H_n(\xi,\xi_i)
K_n(\xi_k\cdot\xi_i)\right|\\
&=&\sum_{k=1}^N \left|\sum_{i=1}^{N}(\lambda_i-\lambda_{i-1})\sum_{j=i}^N H_n(\xi,\xi_j)K_n(\xi_k\cdot\xi_j)\right|
\\
&\le&\sum_{i=1}^{N}(\lambda_i-\lambda_{i-1})\sum_{k=1}^N \left|\sum_{j=i}^N H_n(\xi,\xi_j)K_n(\xi_k\cdot\xi_j)\right|
\\
&\le& \left(\max_{1\le i\le N}\sum_{k=1}^N\left|\sum_{j=i}^N H_n(\xi,\xi_j)
K_n(\xi_k\cdot\xi_j)\right|\right)\left(\sum_{i=1}^{N} \lambda_i-\lambda_{i-1} \right)
\\
&=& \lambda_N \left(\max_{1\le i\le N}\sum_{k=1}^N\left|\sum_{j=i}^N H_n(\xi,\xi_j)
K_n(\xi_k\cdot\xi_j)\right|\right)
\\
&\le& \frac\C{n^2} \left(\max_{1\le i\le N}\sum_{k=1}^N\left|\sum_{j=i}^N H_n(\xi,\xi_j)
K_n(\xi_k\cdot\xi_j)\right|\right).
\end{eqnarray*}
Hence,
(\ref{tesi}) is proved if we show that
\begin{equation}\label{tesi1}
A_i:=\frac 1{n^2}\sum_{k=1}^N\left|\sum_{j=i}^N H_n(\xi,\xi_j)
K_n(\xi_k\cdot\xi_j)\right|\le\C\sqrt{n}
\end{equation}
holds for all $i=1,\ldots,N$, with $\C\ne\C(i,n,N,\xi)$.

We prove (\ref{tesi1}) by induction on $i$ and from now on in this proof we always mean that $\C\ne\C(i,n,N,\xi)$.

For $i=1$, by means of (\ref{inva-LSker}), (\ref{Marci}) and (\ref{Leb-Fou}), we have
\begin{eqnarray*}
A_1 &=& \frac 1{n^2}\sum_{k=1}^N\left|\sum_{j=1}^N H_n(\xi,\xi_j)
K_n(\xi_k\cdot\xi_j)\right|=
\frac 1{n^2}\sum_{k=1}^N\left|K_n(\xi_k\cdot\xi)\right|\\
&\le&  \C\int_{\SS^2}\left|K_n(\y\cdot\xi)\right|d\sigma(\y)\le\C\sqrt{n}.
\end{eqnarray*}
Now, let us assume that $A_i\le\C\sqrt{n}$ and prove $A_{i+1}\le\C\sqrt{n}$.

Indeed, recalling that $\xi\in X_N $, we can apply (\ref{ls-ker1}), and using also (\ref{Marci}) and (\ref{Leb-Fou}), we get
\begin{eqnarray*}
A_{i+1} &=& \frac 1{n^2}\sum_{k=1}^N\left|\sum_{j=i}^N H_n(\xi,\xi_j)
K_n(\xi_k\cdot\xi_j)-H_n(\xi,\xi_i)K_n(\xi_k\cdot\xi_i)\right|\\
&\le& A_i +\frac 1{n^2}\sum_{k=1}^N\left|H_n(\xi,\xi_i)K_n(\xi_k\cdot\xi_i)\right|\le A_i+\frac 1{n^2}\sum_{k=1}^N\left|K_n(\xi_k\cdot\xi_i)\right|\\
&\le& A_i+\C\int_{\SS^2}\left|K_n(\y\cdot\xi_i)\right|d\sigma(\y)\le A_i+\C\sqrt{n}\le\C\sqrt{n}.
\end{eqnarray*}
This proves the statement (\ref{tesi}) in the case that $\xi\in X_N$.

For arbitrary $\xi\in\SS^2$, we reason analogously, but we start applying (\ref{inva-Fou}) to the polynomials $P(\x)=H_n(\x,\xi_k)$, with $k=1,\ldots,N$. In this way, by (\ref{sum-part}) and (\ref{l}) we get
\begin{eqnarray*}
\sum_{k=1}^N|H_n(\xi,\xi_k)|&\le&\sum_{k=1}^N\left|\sum_{i=1}^N\lambda_i  H_n(\xi_i,\xi_k)
K_n(\xi\cdot\xi_i)\right|\\
&\le&
\lambda_N\left(\max_{1\le i\le N}\sum_{k=1}^N\left|\sum_{j=i}^N H_n(\xi_j,\xi_k)K_n(\xi\cdot\xi_j)\right|\right)
\\
&\le& \frac\C{n^2}\left(\max_{1\le i\le N}\sum_{k=1}^N\left|\sum_{j=i}^N H_n(\xi_j,\xi_k)
K_n(\xi\cdot\xi_j)\right|\right).
\end{eqnarray*}
Then, set
\[
B_i:=\frac 1{n^2}\sum_{k=1}^N\left|\sum_{j=i}^N H_n(\xi_j,\xi_k)K_n(\xi\cdot\xi_j)\right|,\qquad i=1,\ldots,N,
\]
we note that similarly to $A_1$ we have
\[
B_1=\frac 1{n^2}\sum_{k=1}^N\left|\sum_{j=1}^N H_n(\xi_j,\xi_k)K_n(\xi\cdot\xi_j)\right|=
\frac 1{n^2} \sum_{k=1}^N\left|K_n(\xi\cdot\xi_k)\right|\le\C\sqrt{n}.
\]
Moreover, recalling (\ref{sup-Darboux})
and taking into account that we have already proved (\ref{tesi}) with $\xi=\xi_i\in X_N$, we get
\[
B_{i+1}\le B_i+\frac 1{n^2}\sum_{k=1}^N\left|H_n(\xi_i,\xi_k)K_n(\xi\cdot\xi_i)\right|\le
B_i+\C\sum_{k=1}^N\left|H_n(\xi_i,\xi_k)\right|\le B_i +\C\sqrt{n}.
\]
Hence the statement (\ref{tesi}) follows by induction in the case $\xi\in\SS^2$ too.
\subsection{Proof of Theorem \ref{th-VPls}}
Let be arbitrarily fixed $0<\theta<1$ and $f$ such that $\|f\|_\infty<\infty$. Set $m=\lfloor \theta n\rfloor$ and
\begin{equation}\label{LSmean-ker}
\tilde v_n^m(\x,\y):=\frac 1{2m+1}\sum_{r=n-m}^{n+m}H_r(\x,\y),\qquad \x,\y\in\SS^2,
\end{equation}
by (\ref{meanLS}) and (\ref{LS-sum}), we have
\begin{equation}\label{LSmean-sum}
\tilde V_n^mf(\x)=\sum_{k=1}^Nf(\xi_k)\tilde v_n^m(\x, \xi_k),\qquad \x\in\SS^2.
\end{equation}
Taking into account that the assumptions of Lemma \ref{lem-quad} are satisfied with $n$ replaced by $2n$, and that
\[
\deg{\left(\tilde V_n^mf\right)}=n+m=\lfloor (1+\theta)n\rfloor=\left\lfloor\frac{1+\theta}2 \ 2n\right\rfloor <2n,
\]
we can apply (\ref{Marci-inf}) to $P=\tilde V_n^mf\in\PP_{\lfloor 2n\rho\rfloor}$ with $\rho=(1+\theta)/2$. In this way,  by (\ref{LSmean-sum}) we deduce that
\begin{eqnarray*}
\|\tilde V_n^mf\|_\infty &\le& \C \max_{1\le i\le N}|\tilde V_n^mf(\xi_i)|
= \C \max_{1\le i\le N} \left|\sum_{k=1}^Nf(\xi_k)\tilde v_n^m(\xi_i, \xi_k)\right|\\
&\le& \C \left(\max_{1\le k\le N}|f(\xi_k)|\right)
\left(\max_{1\le i\le N}\sum_{k=1}^N\left|\tilde v_n^m(\xi_i, \xi_k)\right|\right).
\end{eqnarray*}
Hence, we are going to get the statement  by proving that
\begin{equation}\label{tesi2}
\Sigma:=\sum_{k=1}^N\left|\tilde v_n^m(\xi, \xi_k)\right|\le\C, \qquad \forall \xi\in X_N,
\end{equation}
where from now on in this proof we assume that $\C\ne\C(n,N,\xi)$.

Let $\xi\in X_N$ be arbitrarily fixed. Note that for $r=n-m,\ldots,n+m$, if we use (\ref{inva-Fou}) with $P(\x)=H_r(\xi,\x)$, and then apply (\ref{quad}) with $P(\x)=H_r(\xi,\x)K_r(\xi_k\cdot\x)\in\PP_{2(n+m)}\subset\PP_{4n}$, we obtain
\[
H_r(\xi,\xi_k)=\int_{\SS^2}H_r(\xi,\y)K_r(\xi_k\cdot\y)d\sigma(\y)=\sum_{i=1}^N\lambda_iH_r(\xi,\xi_i)
K_r(\xi_k\cdot\xi_i),\qquad k=1,\ldots,N.
\]
Consequently,  supposing that (\ref{l}) holds, if we use (\ref{sum-part}) (with $a_i=\lambda_i$)  and  (\ref{li}), we get
\begin{eqnarray*}
\Sigma&:=&\sum_{k=1}^N\left|\frac 1{2m+1}\sum_{r=n-m}^{n+m}H_r(\xi, \xi_k)\right|\\
&=& \sum_{k=1}^N\left|\frac 1{2m+1}\sum_{r=n-m}^{n+m}\sum_{i=1}^N\lambda_i H_r(\xi,\xi_i)
K_r(\xi_k\cdot\xi_i)\right|\\
&=& \sum_{k=1}^N\left|\sum_{i=1}^N\lambda_i\left[\frac 1{2m+1}\sum_{r=n-m}^{n+m}H_r(\xi,\xi_i)
K_r(\xi_k\cdot\xi_i)\right]\right|\\
&\le& \sum_{k=1}^N\sum_{i=1}^N(\lambda_i-\lambda_{i-1})\left|\sum_{j=i}^N\left[\frac 1{2m+1}\sum_{r=n-m}^{n+m}H_r(\xi,\xi_j)
K_r(\xi_k\cdot\xi_j)\right]\right|
\\
&\le&\lambda_N\left(\max_{1\le i\le N}\sum_{k=1}^N\left|\sum_{j=i}^N\left[\frac 1{2m+1}\sum_{r=n-m}^{n+m}H_r(\xi,\xi_j)
K_r(\xi_k\cdot\xi_j)\right]\right|\right)
\\
&\le&\frac\C{n^2}\left(\max_{1\le i\le N}\sum_{k=1}^N\left|\sum_{j=i}^N\left[\frac 1{2m+1}\sum_{r=n-m}^{n+m}H_r(\xi,\xi_j)
K_r(\xi_k\cdot\xi_j)\right]\right|\right).
\end{eqnarray*}
Thus, similarly to the previous subsection, we set
\[
A_i:=\frac 1{n^2}\sum_{k=1}^N\left|\sum_{j=i}^N\left[\frac 1{2m+1}\sum_{r=n-m}^{n+m}H_r(\xi,\xi_j)
K_r(\xi_k\cdot\xi_j)\right]\right|,\qquad i=1,\ldots, N,
\]
and by induction on $i$ we are going to prove that $A_i\le\C$.

For $i=1$, by means of (\ref{inva-LSker}),  (\ref{Marci}) and (\ref{LC-VPclas}), we get
\begin{eqnarray*}
A_1 &:=& \frac 1{n^2}\sum_{k=1}^N\left|\frac 1{2m+1}\sum_{r=n-m}^{n+m}\left(\sum_{j=1}^NH_r(\xi,\xi_j)
K_r(\xi_k\cdot\xi_j)\right)\right|\\
&=& \frac 1{n^2}\sum_{k=1}^N\left|\frac 1{2m+1}\sum_{r=n-m}^{n+m}
K_r(\xi_k\cdot\xi)\right|\\
&\le&\C\int_{\SS^2}\left|\frac 1{2m+1}\sum_{r=n-m}^{n+m}
K_r(\y\cdot\xi)\right|d\sigma(\y)\le\C.
\end{eqnarray*}
In order to prove that $A_i\le\C\Longrightarrow A_{i+1}\le \C$, we note that by (\ref{Marci})
\begin{eqnarray*}
A_{i+1} &:=& \frac 1{n^2}\sum_{k=1}^N
\left|\frac 1{2m+1}\sum_{r=n-m}^{n+m}\left(\sum_{j=i}^NH_r(\xi,\xi_j)
K_r(\xi_k\cdot\xi_j)- H_r(\xi,\xi_i)
K_r(\xi_k\cdot\xi_i)\right)\right|\\
&\le & A_i+\frac 1{n^2}\sum_{k=1}^N\left|\frac 1{2m+1}\sum_{r=n-m}^{n+m}
H_r(\xi,\xi_i)K_r(\xi_k\cdot\xi_i)\right|\\
&\le & A_i+ \C\int_{\SS^2}\left|\frac 1{2m+1}\sum_{r=n-m}^{n+m}
H_r(\xi,\xi_i)K_r(\y\cdot\xi_i)\right|d\sigma(\y).
\end{eqnarray*}
On the other hand,
by means of the following summation by part formula
\[
\sum_{r=\nu}^\mu a_rb_r=a_\mu\sum_{r=\nu}^\mu b_r+\sum_{r=\nu}^{\mu-1}(a_r-a_{r+1})\sum_{s=\nu}^rb_r,
\]
for any $i=1,\ldots,N$,
we have
\begin{eqnarray*}
\sum_{r=n-m}^{n+m}
H_r(\xi,\xi_i)K_r(\y\cdot\xi_i)&=& H_{n+m}(\xi,\xi_i)\sum_{r=n-m}^{n+m}
K_r(\y\cdot\xi_i)+\\
&& +\sum_{r=n-m}^{n+m-1}\left[H_r(\xi,\xi_i)-H_{r+1}(\xi,\xi_i)\right]
\sum_{s=n-m}^rK_s(\y\cdot\xi_i),
\end{eqnarray*}
and consequently by (\ref{ls-ker1}) and Lemma \ref{lem-gen} we get
\begin{eqnarray*}
&&\int_{\SS^2}\left|\frac 1{2m+1}\sum_{r=n-m}^{n+m}
H_r(\xi,\xi_i)K_r(\y\cdot\xi_i)\right|d\sigma(\y)\\
&\le & |H_{n+m}(\xi,\xi_i)|\int_{\SS^2}\left|\frac 1{2m+1}\sum_{r=n-m}^{n+m}
K_r(\y\cdot\xi_i)\right|d\sigma(\y)\\
&&+ \sum_{r=n-m}^{n+m-1}\left|H_r(\xi,\xi_i)-H_{r+1}(\xi,\xi_i)\right|
\int_{\SS^2}\left|\frac 1{2m+1}\sum_{s=n-m}^{r}
K_s(\y\cdot\xi_i)\right|d\sigma(\y)\\
&\le&
\left(\sup_{n-n\le r\le n+m} \int_{\SS^2}\left|\frac 1{2m+1}\sum_{s=n-m}^{r}
K_s(\y\cdot\xi_i)\right|d\sigma(\y)\right)\cdot\\
&&\cdot \left(|H_{n+m}(\xi,\xi_i)|+\sum_{r=n-m}^{n+m-1}\left|H_r(\xi,\xi_i)-H_{r+1}(\xi,\xi_i)\right|\right)\\
&\le& \C \left(|H_{n+m}(\xi,\xi_i)|+\sum_{r=n-m}^{n+m-1}\left|H_r(\xi,\xi_i)-H_{r+1}(\xi,\xi_i)\right|\right)\\
&=& \C \left|\sum_{k=1}^{(n+m+1)^2}I_k(\xi)I_k(\xi_i)\right|+
\C\sum_{r=n-m}^{n+m-1}\left|\sum_{k=(r+1)^2+1}^{(r+2)^2}I_k(\xi)I_k(\xi_i)\right|\\
&\le& \C\sum_{k=1}^{(n+m+1)^2}\left|I_k(\xi)I_k(\xi_i)\right|.
\end{eqnarray*}
Summing up, we have got
\begin{eqnarray*}
A_{i+1}&\le& A_i+ \C\int_{\SS^2}\left|\frac 1{2m+1}\sum_{r=n-m}^{n+m}
H_r(\xi,\xi_i)K_r(\y\cdot\xi_i)\right|d\sigma(\y)\\
&\le& A_i + \C\sum_{k=1}^{(n+m+1)^2}\left|I_k(\xi)I_k(\xi_i)\right|,
\end{eqnarray*}
where by our assumptions and by Lemma \ref{lem-quad} (with $n$ replaced by $2n$) it is $(n+m+1)^2<(2n+1)^2<N$. Consequently  (\ref{Q}) holds with $n$ replaced by $n+m$ and for all $\xi\in X_N$, we have
\[
\sum_{k=1}^{(n+m+1)^2}\left|I_k(\xi)I_k(\xi_i)\right|\le \left(\sum_{k=1}^{(n+m+1)^2}|I_k(\xi)|^2\right)^\frac 12
\left(\sum_{k=1}^{(n+m+1)^2}|I_k(\xi_i)|^2\right)^\frac 12\le 1.
\]
Hence we deduce that $A_{i+1}\le A_i+\C$, with $\C\ne\C(i,n,N,\xi)$.

Consequently 
$ \limsup_{n\rightarrow \infty} A_{i+1}=+\infty$ implies $\limsup_{n\rightarrow \infty} A_{i}=+\infty$, i.e., the statement follows.
\Proofend
\section{Conclusions}
We considered the polynomial approximation of a function $f$ on the unit Euclidean sphere $\SS^2\subset\RR^3$, by using samples of $f$ at a discrete point set $X_N=\{\xi_1,\ldots,\xi_N\}\subset\SS^2$.

Under the assumption that the points in $X_N$ satisfy the Marcinkiewicz type inequality (\ref{Marci-1}) and they are nodes of a  positive weighted quadrature rule of suitable degree of exactness, we proved that we can use the least squares polynomial approximant $\tilde S_nf$ in (\ref{LS-min}) and the mean $\tilde V_n^mf$ of least squares approximants given in (\ref{meanLS}) as an alternative to the hyperinterpolation polynomial $L_nf$ in (\ref{hyper}) and  the filtered hyperinterpolation polynomial $V_n^mf$ in (\ref{filtered}) respectively.

Indeed by Theorems \ref{th-ls} and \ref{th-VPls}, we showed that they have a comparable approximation degree w.r.t.\ uniform norms, that is, denoted by $E_n(f)$  the error of best uniform approximation of $f$ by spherical polynomials of degree at most $n$,  similarly to hyperinterpolation, we have
\[
E_{n}(f)\le\|f-\tilde\S_nf\|_\infty\le  \C \sqrt{n} E_{n}(f),\qquad \C\ne\C(n,f),
\]
and analogously to filtered hyperinterpolation, for any $0<\theta<1$ and $m=\lfloor \theta n\rfloor$, we get
\[
E_{n+m}(f) \le\|f-\tilde V_n^mf\|_\infty\le  \C E_{n-m}(f),\qquad \C\ne\C(n,f).
\]
These theoretical results were illustrated by numerical experiments.

In the case of functions almost everywhere smooth, apart from some isolated points of singularity, it is known that the Gibbs phenomenon occurs by using least squares as well as hyperinterpolation polynomials. By a numerical experiment we showed that this phenomenon can be strongly reduced if we consider the mean $\tilde V_n^mf$. Keeping $n$  unchanged (which is strictly related to the number $N$ of data) we appropriately modulate the range of action $m$ of the mean by suitably varying the parameter $\theta$ from the limiting value $\theta=0$, which corresponds to simple least squares approximation (suggested for very smooth functions) to the limiting value $\theta =1$, which in practice corresponds to the Fej\'er mean of least squares polynomials.
\vspace{.5cm}\newline
{\bf Acknowledgments }
The authors would like to thank the anonymous referee for his valuable comments and suggestions to improve the quality of the paper, and Ed Saff for providing us reference \cite{HardMichSaff2016}.
The research of the first author was partially supported by GNCS--INDAM, and of the second author by
  the Research Council KU Leuven,
C1-project (Numerical Linear Algebra and Polynomial Computations),
and by
the Fund for Scientific Research--Flanders (Belgium),
``SeLMA'' - EOS reference number: 30468160.
\bibliographystyle{abbrv}
\bibliography{longstrings,TOTAL}

\end{document}